\documentclass[12pt,reqno,a4paper]{amsart}
\usepackage{amsmath,amssymb,amsfonts,amsthm,amscd}
\usepackage[mathscr]{eucal}
\usepackage{hyperref}

\author{Theodore Voronov}

\address{Department of Mathematics, University of Manchester Institute of Science and Technology
(UMIST), United Kingdom}

\email{theodore.voronov@umist.ac.uk}

\title[On inverse problem of calculus of
variations]{An alternative form of the Helmholtz  criterion in the
inverse problem of calculus of variations}


\dedicatory{To Alan Weinstein on the occasion of his 60th birthday}

\newtheorem{thm}{Theorem}
\theoremstyle{definition}
\theoremstyle{remark}
\newtheorem{Rem}{Remark}
\newtheorem*{exe}{Example}
\newtheorem*{ack}{Acknowledgement}

\def\co{\colon\thinspace}
\newcommand{\mathcall}{\EuScript}

\newcommand{\rhs}{r.h.s.}

\renewcommand{\geq}{\geqslant}

\newcommand{\E}{{\mathcall{E}}}

\newcommand{\der}[2]{{\frac{\partial {#1}}{\partial {#2}}}}
\newcommand{\lder}[2]{{\partial {#1}/\partial {#2}}}
\newcommand{\var}[2]{{\frac{\delta {#1}}{\delta {#2}}}}
\newcommand{\lvar}[2]{{\delta {#1}/\delta {#2}}}

\newcommand{\RR}{\mathbb R}

\def\s{\sigma}

\def\O{\Omega}

\def\o{\omega}

\newcommand{\h}{{\eta}}

\def\d{\delta}

\newcommand{\du}{{\delta u}}

\begin{document}

\begin{abstract}
We give a necessary and sufficient condition for the existence of a
local solution of the inverse problem of calculus of variations in
terms of the identical vanishing of the variation of a functional on
an extended space (with the number of independent variables increased
by one), and explain its relation with the classical Helmholtz
criterion using the de Rham complex on an infinite-dimensional space
of fields.
\end{abstract}

\maketitle

\section{Introduction}

This paper deals with the question, under which conditions given
functions can be the variational derivatives of some functional.
(The precise setup is described in \S\ref{secmain}.) This is a
classical question, and there is a classical answer to it given by
the so-called Helmholtz (or Helmholtz--Volterra) criterion. One
can see~\cite{olver} for an exposition and a historical review.
The Helmholtz condition is, in fact, nothing but the closedness of
a $1$-form of a particular appearance on an infinite-dimensional
``space of fields'' on which functionals under consideration are
defined, though this may be hidden in some expositions. Hence the
Helmholtz criterion can be viewed as a version of an
infinite-dimensional Poincar\'e lemma for $1$-forms, and proved
accordingly.

In the ordinary finite-dimensional case we know that a $1$-form is
closed if its integrals over   paths do  not change under small
perturbation of a path fixing the boundary. This is also true for
$k$-forms and $k$-paths. Pushed to the limit, this idea allows to
express the usual de Rham differential entirely in variational
terms. (For supermanifolds this is, actually, the only way one can
construct an adequate analog of the Cartan--de Rham complex,
see~\cite{tv:compl, tv:git, tv:cartan1},  also~\cite{hov:gayduk},
\cite{ass:bar}, \cite{hov:khn2}.) Moreover, the de Rham complex is
embedded in a larger complex consisting of \textit{all}
Lagrangians of $k$-paths on a manifold~\cite{tv:lag}.

In this note we suggest an alternative criterion of the existence
of a local solution of the inverse problem of calculus of
variations in terms of the identical vanishing of the variation of
a functional on an extended space  (Theorem~\ref{thmcrit}).  This
criterion is equivalent to the classical Helmholtz one. One
advantage of the suggested criterion is the simplicity of
application, due to the fact that, roughly speaking, it requires
calculating the differential of a \textit{function}, while the
classical Helmholtz test requires calculating  the differential of
a \textit{$1$-form}. Theorem~\ref{thmcrit} is an immediate
consequence of the existence of a natural chain transformation
diminishing degrees by one and peculiar for forms on field spaces
(Theorem~\ref{thmdiagram}). In a way, this is an ``integrated''
Cartan homotopy formula, with no analog on finite-dimensional
manifolds. (On the other hand, the expression of the usual de Rham
differential via the variation of functionals also follows from
Theorem~\ref{thmdiagram}.)

We give a very brief outline of the appropriate de Rham complex in
the Appendix. Methods of supermathematics are helpful, as usual.
Clearly, forms on infinite-dimensional spaces like spaces of fields
are no novel for physicists. Some mathematicians, on the other hand,
seem to refrain from using them as legitimate objects. However, as we
show, using them explicitly is very convenient and is entirely
rigorous. I would like to note that a formalization based on the
so-called functional forms on an infinite jet space
(see~\cite{olver}, \cite{vinog:kollektiv}), meant to replace the
supposedly non-rigorous differential forms on field spaces, has a
ineradicable defect of modelling only a subspace of the ``diagonal
forms'' (see below). This results in the impossibility to multiply
such objects and in other counterintuitive features. Working with
forms as such is much better!

\section{Main statements}\label{secmain}

Let $E\to M$ be a smooth fiber bundle over a smooth manifold $M$.
The inverse problem of calculus of variations is to determine
whether functions $A_i(x,u,u',u'',\ldots)$  can be the left hand
sides of the Euler--Lagrange equations for a Lagrangian
$L=L(x,u,u',\ldots)$, and to find $L$ if it exists. We denote by
$x=(x^a)$ local coordinates on $M$, by $u=(u^i)$ local coordinates
in the fiber of $E\to M$; then $(x,u,u', u'',
\ldots)=(x^a,u^i,u^i_a,u^i_{ab},\ldots)$ are the coordinates of
jets of sections of $E\to M$. In the sequel we use the notation
such as $D_a$ for total derivatives. The classical Helmholtz
criterion (which we will recall below) states that, locally, a
necessary and sufficient condition of the equality
$A_i=\lvar{S}{u^i}$ for some $S=\int L\, d^nx$ is the
self-adjointness of a certain differential operator constructed
from the functions $A_i$.  We can give an alternative criterion:

\begin{thm}\label{thmcrit}
Functions $A_i(x,u,u',u'',\ldots\,)$ are variational derivatives of
some functional $\int L\,d^nx$, with a Lagrangian $L$ defined locally
in $u^i$,  if and only if,  for the  functional
$\int \dot u^i A_i(x,u,u',u'',\ldots)\,d^nxdt$,
defined on sections of the induced bundle $E\times \RR \to M\times
\RR$, its variation  vanishes identically for all arguments. (Dot
stands for the derivative in $t$.)
\end{thm}

\begin{proof}
Necessity. Let $A_i$ be the variational derivatives of some
functional $S=\int L\,d^nx$. Then we have $L'=\dot u^iA_i=\dot
u^i\,\lvar{S}{u^i}=D_t L +D_af^a$ for some functions
$f^a(x,u,u',\ldots)$, by the definition of variational derivative.
Hence $\int L'\,d^nxdt$ is the integral of a total divergence, and
its variation identically vanishes. Sufficiency. First we shall write
down our condition
\begin{equation}\label{eqcond1}
\delta \int \dot u^i A_i\,d^nxdt=0
\end{equation}
explicitly. By expanding and integrating by parts w.r.t. $x^a, t$ one
can arrive at the following formula:
\begin{multline*}
\delta \int \dot u^i A_i\,d^nxdt= 
\int d^nxdt \,\delta u^i(x,t) \sum\left(- \der{A_i}{u^j_{a_1\ldots
a_k}}\,D_{a_1}\ldots D_{a_k} + \right.
\\
\left.  (-1)^k D_{a_1}\ldots D_{a_k} \circ \der{A_j}{u^i_{a_1\ldots
a_k}}\right) \dot u^j.
\end{multline*}
The symbol $\circ$ stands for the composition of operators (here we
first apply the multiplication by a function and then the
differentiation). Since both $\delta u^i(x,t)$ and $\dot u^i(x,t)$
can be arbitrary functions, we obtain the condition in the form
\begin{equation} \label{eqcond2}
\sum_{k\geq 0}\left(\der{A_i}{u^j_{a_1\ldots a_k}}\,D_{a_1}\ldots
D_{a_k} - \right.
\left.  (-1)^k D_{a_1}\ldots D_{a_k} \circ \der{A_j}{u^i_{a_1\ldots
a_k}}\right)=0.
\end{equation}
Basically, there is nothing to prove now, since~\eqref{eqcond2} is
exactly the classical Helmholtz condition (see below). For the sake
of completeness we shall supply a standard type argument. In a
star-shaped domain in $u^i$, one has $A_i=\lvar{S}{u^i}$ for $S=\int
L \,d^nx$, where
\begin{equation*}
L(x,u,u',u'',\ldots)=\int_0^1 dt\, u^i A_i(x,tu,tu',tu'',\ldots).
\end{equation*}
Indeed, by taking the variation and integrating by parts in $x^a$,
\begin{multline*}
    \delta S=
     \int d^nx \int_0^1\!dt\,\du^i
    \Biggl(A_i(x,tu,\ldots) \, + \\
      t \sum (-1)^k D_{a_1}\ldots D_{a_k}
    \left(\der{A_j}{u^i_{a_1\ldots a_k}}\,(x,tu,\ldots)\,
    u^j\right)\Biggl),
\end{multline*}
{which by~\eqref{eqcond2} equals}
\begin{multline*}
    \int d^nx \int_0^1\!dt\,\du^i
    \Biggl(A_i(x,tu,\ldots)  +
      t \sum
     \der{A_i}{u^j_{a_1\ldots a_k}}\,(x,tu,\ldots)\, u^j_{a_1\ldots a_k}
     \Biggl)=\\
    \int d^nx \,\du^i  \int_0^1\!dt\,
    \frac{d}{dt}\Bigl(tA_i(x,tu,\ldots)\Bigr)=
    \int d^nx  \,\du^i
      A_i(x,u,\ldots)
\end{multline*}
(we have suppressed the arguments in $u^i(x), \du^i(x)$, etc.).
\end{proof}

To practically apply this theorem one has to calculate the
Euler--Lagrange expression for $L'=\dot u^i A_i$, which is linear in
the  derivatives involving ``time'': $\dot u^i$, $\dot u^i_a$, etc.,
and set the respective coefficients to zero. This will give the
equations for $A_i$.

\begin{exe}
For a second-order function $A=A(x,u,u',u'')$ in the  scalar case the
only non-trivial equation will be
\begin{equation*}
    \der{A}{u_a}=D_b\left(\der{A}{u_{ba}}\right).
\end{equation*}
\end{exe}

The necessity statement in Theorem~\ref{thmcrit} is just the equality
$d^2=0$ in the ``complex of Lagrangians'' introduced in~\cite{tv:lag}
(see also~\cite{tv:varcomplexes}).

The sufficiency is more delicate.

Recall the classical \textit{Helmholtz criterion} (see, e.g.,
\cite{olver}): for functions $A_i(x,u,\ldots)$ to be variational
derivatives, a necessary and sufficient condition is the formal
self-adjointness of the matrix differential operator $L_{ij}$
associated with $A_i$ by the formula $L_{ij}=\sum
\lder{A_j}{u^i_{a_1\ldots a_k}}\,D_{a_1}\ldots D_{a_k}$, i.e.,
\begin{equation}\label{eqcond3}
    L_{ij}=L^*_{ji}.
\end{equation}
Written explicitly, \eqref{eqcond3} is exactly
equation~\eqref{eqcond2}. (Olver~\cite{olver} deduces the local
sufficiency of~\eqref{eqcond3} from properties of the variational
complex on jet space. A direct prove is included above for
completeness.)

The main observation is that after calculating the variation, the
condition~\eqref{eqcond1} reduces to the same
equation~\eqref{eqcond2}, hence~\eqref{eqcond1} and
~\eqref{eqcond3} are equivalent.

A deeper explanation can be given as follows.

The self-adjointness condition~(\ref{eqcond3}), (\ref{eqcond2}) is
nothing but the closedness of a $1$-form  of a special appearance
on the infinite-dimensional ``space of fields'' $u=\{u^i(x)\}$,
i.e., the space of sections of $E$.   The Helmholtz criterion can
be seen as a special case of the Poincar\'e lemma in this
infinite-dimensional situation.

Indeed, for the  differential  of a $1$-form $A=\int\!
d^nx\,\du^i(x)\,A_i(x;[u])$ we have
\begin{equation*}
\d A=\frac{1}{2}\int\!\! d^nx\,d^ny
\;\du^i(x)\,\du^j(y)\left(\var{A_j(y;[u])}{u^i(x)}-\var{A_i(x;[u])}{u^j(y)}\right).
\end{equation*}
Hence $\d A=0$ means
\begin{equation}\label{eqclosed}
\var{A_j(y;[u])}{u^i(x)}-\var{A_i(x;[u])}{u^j(y)}=0.
\end{equation}
(See Appendix for a description of the corresponding de Rham
complex.)

Now, functions $A_i(x,u,u',\ldots)$ can be viewed as the
coefficients of a ``diagonal'' $1$-form $A=\int \! d^nx \;
\du^i(x)\,A_i(x,u,u',\ldots)$, see Appendix. It is easy to check
by expanding the variational derivatives in~\eqref{eqclosed} that
for such kind of $1$-forms the condition $\d A=0$
gives~\eqref{eqcond2}. The differential operators appearing
in~\eqref{eqcond2} arise from derivatives of $\d$-functions. The
proof of the local existence of $S$ simply  follows the usual
proof of the Poincar\'e lemma.

The condition~\eqref{eqcond1}, on the other hand, is the closedness
of a  $0$-form. We will see that the fact that ~\eqref{eqcond1} and
\eqref{eqcond3} give the same thing follows from the commutativity of
the differentials in the de Rham complexes on field  spaces  with a
natural homomorphism $K$ defined below.

Let $\E(M)$ stand for the space of sections of the bundle $E\to M$,
and $\E(M\times \RR)$ for the  space of sections of the induced
bundle over $M\times \RR$. Let $\O(\E(M))$ and $\O(\E(M\times \RR))$
denote the corresponding algebras of forms. There is a natural odd
map $K\co \O(\E(M))\to \O(\E(M\times \RR))$ that lowers the degree by
one:
\begin{equation*}
K\o:=\int d^nx\, dt \;\dot u^i(x,t)\,\var{\o}{\du^i(x)}.
\end{equation*}
At the {\rhs} we treat sections on $M\times \RR$ as families of
sections on $M$. $K$ is the composition of the map $I=\int dt\co
\O(\E(M))\to \O(\E(M\times \RR))$ that sends every functional to its
integral over $t$ and the interior product with the canonical vector
field $\dot u$ on $\E(M\times \RR)$. Notice that $K$ is monomorphic
on forms of degree $>0$.

\begin{thm}\label{thmdiagram}
The following diagram commutes:
\begin{equation*}\label{eqdiagram}
    \begin{CD}
    \O^k(\E(M))@>{K}>> \O^{k-1}(\E(M\times \RR)) \\
    @V{\d}VV                  @VV{\d}V \\
    \O^{k+1}(\E(M))@>{K}>> \O^{k}(\E(M\times \RR))
    \end{CD}
\end{equation*}

\end{thm}
\begin{proof}
Consider the tautological family of maps $f_t\co \E(M\times \RR) \to
\E(M)$ that sends every section over $M\times \RR$ to itself
considered as a family of sections. Then $\dot u$ can be
alternatively viewed as the time-dependent velocity vector field for
$f_t$, so that $i_{\dot u}\co \O(\E(M))\to \O(\E(M\times \RR))$. The
homotopy formula gives
\begin{equation*}
\left(\d\, i_{\dot u} +i_{\dot u}\, \d\right)\o
=\frac{d}{dt}\,f^*_t\o
\end{equation*}
for any form $\o$ on $\E(M)$. Hence by integrating w.r.t. $t$ we
get
\begin{equation}\label{eqcommut}
\d K+ K\d=0,
\end{equation}
as claimed.  The {\rhs} of~\eqref{eqcommut} is zero as the
integral of a total derivative. The commutator has  a plus sign
since both $\d$ and $K$ are odd.
\end{proof}

In particular, for $k=1$, we get the relation between the
differentials of a $1$-form $A$ and the $0$-form $KA$: since $K$
is monomorphic, $\d A=0$ if and only if $K(\d A)=-\d (KA)=0$. This
is precisely the relation between the classical Helmholtz
criterion and our Theorem~\ref{thmcrit}.

\begin{Rem}
Suppose $E$ is a bundle over a point, i.e., $M=\{*\}$. Then
$\O(\E(M))=\O(E)$, and Theorem~\ref{thmdiagram} relates forms on a
finite-dimensional manifold $E$ with forms on the space of paths.
(Without harm,  $\RR$ can be replaced above by $I=[0,1]$.)
Iterating, we get the commutative diagram
\begin{equation*}\label{eqdiagram2}
    \begin{CD}
    \O^k(E)@>{K^k}>> \O^{0}(\{I^{k}\to E\})\\
    @V{d}VV                  @VV{K\circ \d}V \\
    \O^{k+1}(E)@>{K^{k+1}}>> \O^{0}(\{I^{k+1}\to E\})
    \end{CD}
\end{equation*}
in which we recognize the expression of the differential of
$k$-forms on a manifold $E$ via the variation of functionals of
$k$-paths.
\end{Rem}

\begin{Rem}
With obvious changes, e.g., $d^nx$  replaced by the Berezin volume
element, all  statements remain true for supermanifolds.
\end{Rem}

\appendix
\section*{Appendix. The de Rham complex on a space of fields}

We can consider on the space of sections of $E\to M$ functionals of
the following form: $F[u]=\int f \,d^nx_1\ldots d^nx_k$ and their
sums, where the (formal) integration is over $M\times \ldots \times
M$ and the integrand $f$ is allowed to depend on a finite number of
derivatives of $u^i$ at the points $x_1, \ldots, x_k$, as well as,
possibly, at some other points. Particular examples are the classical
``local functionals'' $S[u]=\int L(x,u(x), \ldots) \,d^nx$ and
``point functionals'' such as $u^i(x)$ or $u^i_a(x)$. One should also
allow in $f$   products of $\delta$-functions and their derivatives
taken at distinct points of $M$; otherwise $f$ should be smooth. Such
functionals make an algebra closed under taking variational
derivatives, which act as derivations. Variational derivatives reduce
the ``integrality'' of a functional (the number of integrations minus
the number of $\delta$-functions involved) by one. Without loss of
generality, integrands can be considered symmetric w.r.t. the
arguments corresponding to the integration points. Then the
variational derivative of any such functional will be given by the
usual Euler--Lagrange expression w.r.t. one point with the
integration remaining over the other points. (Functionals of this
kind have been considered by physicists, see,
e.g.,~\cite{dewitt:dynam}.)

Now, a vector at a ``point'' $\{u^i(x)\}$ is, of course, a section of
$u^*T^{\text{vert}}E$ (the pull-back of the vertical tangent bundle).
A vector field can be formally written as
\begin{equation*}
\h=\int d^nx \,\h^i(x;[u])\,\var{}{u^i(x)},
\end{equation*}
which is a functional of $u$ taking values in $u^*T^{\text{vert}}E$.

The simplest way to define \textit{forms} on this
infinite-dimensional space is to consider functionals of pairs of
fields $(u, \delta u)$, where $\delta u$ is odd. (For
supermanifolds, $\du$ should have parity opposite to that of $u$.)
More precisely, $\delta u$ is a section of $u^*\Pi
T^{\text{vert}}E$, where $\Pi$ is the parity reversion functor.
Forms make an algebra. Analogs of usual formulae hold, viz.,
\begin{equation}\label{eqextdelta}
\delta=\int d^nx \,\delta u^i(x)\,\var{}{u^i(x)}
\end{equation}
for the differential (the ``exterior  variation''),  and
\begin{equation}\label{eqintprod}
i_{\h}=\int d^nx \,\h^i(x;[u])\,\var{}{\delta u^i(x)}
\end{equation}
for the interior product with a vector field $\h$, where at the
{\rhs} of ~\eqref{eqintprod} stands the variational derivative w.r.t.
the odd field $\delta u^i(x)$. One can easily see that the homotopy
formula and the Poincar\'e lemma (with the usual proof) hold as on
ordinary manifolds; in particular, as a form $\s$ such that $\d\s=\o$
if $\d\o=0$ in a star-shaped domain one can take
\begin{equation*}
\s[u,\du]=\int_0^1 dt \int d^nx\,
u^i(x)\,\var{\o}{\du^i(x)}\,[tu,t\du],
\end{equation*}
similarly to the ordinary case.

\begin{Rem}\label{remforms}
``Functional forms'' on the infinite jet space $J^{\infty}(E)$,
considered in \cite{olver}, can be interpreted as corresponding to
a special subclass of all forms on the space of sections.  We call
a $k$-form $\o =\o[u,\du]$ on the space of sections,
\begin{equation}\label{eqform}
\o =\int\!\! d^nx_1\ldots d^nx_k \;\du^{i_1}(x_1)\ldots
\du^{i_k}(x_k) \,\o_{i_1\ldots i_k}(x_1,\ldots,x_k;[u]),
\end{equation}
\textit{diagonal}  if  all coefficients $\o_{i_1\ldots
i_k}(x_1,\ldots,x_k;[u])$ are supported at the diagonal
$x_1=\ldots=x_k=x$ and they are point functionals of $u$ at the
same point $x$. Hence a diagonal form can be re-written
(non-canonically) with only one integration as
\begin{equation*}
\o=\int\!\! d^nx  \;\du^{i_1}_{A_1}(x)\ldots \du^{i_k}_{A_k}(x)
\,\o_{i_1\ldots i_k}^{A_1\ldots A_k}(x,u(x),u'(x),\ldots )
\end{equation*}
where $A_1,\ldots,A_k$ are multi-indices. Thus it can be
identified with a ``functional form'' on  jet space as defined
in~\cite{olver}. ``Functional forms'' appear in the jet-theoretic
approach~\cite{olver}, \cite{vinog:kollektiv} as elements of the
term $E_1$ of the spectral sequence of the bicomplex
$\O^{**}(J^{\infty}(E))$. Effective restriction by diagonal forms,
as well as by similar diagonal multivector fields appearing under
the guise of ``functional multivectors''~\cite{olver},  is a
fundamental  limitation of this approach. The subspace of diagonal
forms is closed under the differential. However, to be able to
multiply forms or consider duality with multivector fields, one
has to work with arbitrary forms.  Even when  results concern only
the diagonal forms, using arbitrary forms gives a clearer picture
and busts intuition.
\end{Rem}

\begin{ack}
I  thank V.~G.~Kac, whose question prompted me to refresh these
topics, and H.~M.~Khudaverdian for discussion. The results of this
note were reported at seminars at the Steklov Mathematical
Institute and the Northeastern University. I am grateful to the
organizers and participants of these seminars, especially to
A.~G.~Sergeev and M.~A.~Shubin. A preliminary version of the text
was contributed to an informal volume dedicated to Alan
Weinstein's 60th birthday.
\end{ack}

\bibliographystyle{hplain}

\def\cprime{$'$} \def\cprime{$'$}

\end{document}